\documentclass{jnmp}

\usepackage{amsmath}
\JNMPnumberwithin{equation}{section}

\newtheorem{theorem}{Theorem}
\newtheorem{proposition}{Proposition}

\theoremstyle{definition}
\newtheorem{remark}{Remark}

\begin{document}

\renewcommand{\evenhead}{S~Bouarroudj}
\renewcommand{\oddhead}{Cohomology of Groups of Diffeomorphisms}

\thispagestyle{empty}

\FirstPageHead{9}{4}{2002}{\pageref{Bouarroudj-firstpage}--
\pageref{Bouarroudj-lastpage}}{Article}

\copyrightnote{2002}{S~Bouarroudj}

\Name{Cohomology of Groups of Diffeomorphisms Related to the Modules of
Differential Operators on a Smooth Manifold}
\label{Bouarroudj-firstpage}

\Author{Sofiane BOUARROUDJ}

\Address{Department of Mathematics, Keio University, Faculty of Science
\& Technology \\
3-14-1, Hiyoshi, Kohoku-ku, Yokohama, 223-8522, Japan\\
E-mail: sofbou@math.keio.ac.jp}

\Date{Received December 21, 2001; Revised April 11, 2002;
Accepted April 17, 2002}

\begin{abstract}
\noindent
Let $M$ be a manifold and $T^*M$ be the cotangent bundle. We introduce a
1-cocycle on the group of diffeomorphisms of $M$ with values in the space of
linear differential operators acting on $C^{\infty} (T^*M).$ When $M$ is
the $n$-dimensional sphere, $S^n$, we use this 1-cocycle to compute the
first-cohomology group of the group of diffeomorphisms of $S^n$, with
coefficients in the space of linear differential operators acting on
contravariant tensor fields.
\end{abstract}

\section{Introduction}
The origin of our investigation is the relationship between the space of
linear differential operators acting on densities and the corresponding space
of symbols, both viewed as modules over the group of diffeomorphisms and the
Lie algebra of smooth vector fields. These two spaces are obviously
isomorphic as vector spaces; however, these spaces are not isomorphic as
modules (cf.~\cite{b,boh,do,lo,lo2}). More precisely, the module of linear
differential operators is a non-trivial deformation of the module of symbols
in the sense of Neijinhuis--Richardson's theory of deformation
(see~\cite{nr}). This theory of deformation can be summarized briefly as
follows. Let a Lie algebra $\mathrm g$ and a ${\mathrm g}$-module $V$
be given. Deformation of the ${\mathrm g}$-module $V$ means we extend the
action of the Lie algebra ${\mathrm g }$ on the formal power series
$V[[t]],$ where $t$ is a parameter. The problem of deformation is related to
the cohomology groups ${\mathrm H}^1({\mathrm g}; \mathrm{End}  (V))$ and
${\mathrm H}^2({\mathrm g}; \mathrm{End}  (V)).$ The first cohomology group
classifies all {\it infinitesimal} deformation of the module $V$, while the
second cohomology group measures the obstruction to extend this infinitesimal
deformation to a formal deformation.

In this paper, we deal with the space of symbols $\mathcal S$; they are
functions on the cotangent bundle which are polynomial on the fibers. This
space is naturally a module over the Lie algebra of smooth vector fields
$\mathrm{Vect}(M)$ and the group of diffeomorphisms $\mathrm{Diff}(M).$ 
Since the module of linear differential operators is a deformation of the 
space of symbols, it is interesting to exhibit 1-cocycles generating the 
infinitesimal deformation. 
For the Lie algebra $\mathrm{Vect}(M)$, Lecomte and Ovsienko~\cite{lo2}
have computed the cohomology group ${\mathrm H}^1(\mathrm{Vect} (M);
{\mathcal D} ({\mathcal S})),$
where ${\mathcal D}({\mathcal S})$ is the space of linear differential 
operators acting on ${\mathcal S}.$ In this paper we will study the 
counterpart cohomology of the group of diffeomorphisms
\begin{equation}
\label{groupe}
{\mathrm H}^1(\mathrm{Diff}(M); {\mathcal D} ({\mathcal S})).
\end{equation}
First, we introduce a non-trivial 1-cocycle on the group of diffeomorphisms
$\mathrm{Diff}(M)$ with values in ${\mathcal D} (C^{\infty}(T^* M)).$ The 
construction of this 1-cocycle was inspired from the Lecomte--Ovsienko's 
result \cite{lo2} where they show that one of the generators of the 
cohomology group 
${\mathrm H}^1(\mathrm{Vect} (M); {\mathcal D} ({\mathcal S}))$ is given by 
the so-called Vey cocycle (see, e.g., \cite{ro})-written with a (lifted) 
connection on the cotangent bundle $T^*M.$ We will explain in the 
Section~\ref{las} the relation between our 1-cocycle and the Vey cocycle. At 
last, we use this 1-cocycle to compute the cohomology group~(\ref{groupe}) 
for the case when the manifold is the $n$-dimensional sphere.

\section{Introducing a 1-cocycle on $\boldsymbol{\mathrm{Diff}(M)}$}
\label{basic}

Let $M$ be a manifold of dimension $n.$ Fix a symmetric affine connection on
it. Let us recall here the natural way to lift a connection to the cotangent
bundle $T^*M$ (see \cite{yano} for more details).

\subsection{Lift of a connection}

Denote by $\Gamma_{ij}^k$, for $i,j,k=1,\ldots,n$ the Christoffel symbols of
the connection on $M.$ There exists a symmetric affine connection on $T^*M$
whose Christoffel symbols $\tilde \Gamma_{\mathbf{ij}}^{\mathbf k}$, where
$\mathbf{i},\mathbf{j}, \mathbf{k}=1, \ldots,2n,$ are as follows: denote
the superscript $\mathbf i$ by $i$ if it varies from $1$ to $n,$ and by 
$\bar{i}$ if it varies from $n+1$ to $2n.$ In local coordinates 
$(x^i, \xi_i)$ on $T^*M,$ the Christoffel symbols of the lifted connection 
are as follows:
\begin{gather}
\tilde \Gamma_{ij}^k =\Gamma_{ij}^k , \qquad  \tilde{\Gamma}_{i{\bar j}}^{k}=0,
\quad \tilde {\Gamma}_{{\bar i}j}^{k}=0, \qquad \tilde{\Gamma}_{{\bar i}
{\bar j}}^{ k}=0,\nonumber\\
\tilde \Gamma_{ij}^{\bar k}= \xi_a (\partial_k \Gamma_{ij}^a-\partial_i
\Gamma_{jk}^a-\partial_j  \Gamma_{ik}^{a}+ 2\Gamma_{kt}^a \Gamma_{ij}^{t}),
\nonumber\\
\tilde \Gamma_{i\bar j}^{\bar k} =-\Gamma_{ik}^j , \qquad
\tilde{\Gamma}_{{\bar i}j}^{\bar k}=-\Gamma_{kj}^i, \qquad
\Gamma_{{\bar i}\bar j}^{\bar k}=0. \label{lift}
\end{gather}
Denote by $\tilde \nabla$ the covariant derivative associated with
the connection $\tilde \Gamma.$

\begin{remark}  The connection $\tilde \Gamma$ is actually the Levi--Civita
connection of some metric on $T^*M$ (see \cite{yano}).
\end{remark}

\subsection{Main definition}

First we shall construct canonically a tensor field on $T^*M.$

It is well known that the difference between two connections is a well-defined
tensor field of type $(2,1).$ It follows from this fact that the following
expression
\begin{equation}
\label{masri}
C (F) := F^* (\tilde  \Gamma)- \tilde \Gamma,
\end{equation}
where $F$ is a lift to $T^*M$ of a diffeomorphism $f\in \mathrm{Diff}(M)$
and $\tilde \Gamma$ is the connection above, is a well-defined tensor field
on $T^*M$ of type $(2,1).$ It is easy to see that the map
\[
F\mapsto C\left(F^{-1}\right)
\]
 defines a non-trivial 1-cocycle on the group of diffeomorphisms of $T^*M$
with values in the space of tensor fields of type $(2,1).$

Denote by $\mathrm{g}^{\mathbf{ij}}$ the standard bivector, dual to the 
standard symplectic structure on $T^*M.$

For all $f\in \mathrm{Diff}(M),$ denote by $\tilde f $ its symplectic lift 
to $T^*M.$

\resetfootnoterule

The main tool of this paper is the third-order linear differential
operator ${\mathcal L}(f)$ on $C^{\infty}(T^* M)$ defined by\footnote{Here and
below we use the convention of summation over repeated indices.}
\begin{gather}
{\mathcal L}(f):=Sym_{\mathbf j,\mathbf i,\mathbf k}
\left (C_{\mathbf m \mathbf l}^{\mathbf j}(\tilde f)
\!\cdot\! \mathrm{g}^{\mathbf{im} }\! \cdot \mathrm{g}^{ \mathbf{kl} }
\right ) \tilde  \nabla_{\mathbf i}  \tilde
\nabla_{\mathbf j} \tilde \nabla_{\mathbf k}\nonumber\\
\phantom{{\mathcal L}(f):=}{}-\frac{3}{2}
Sym_{ \mathbf n,\mathbf m,\mathbf i}
\left (
C_{\mathbf l \mathbf k}^{ \mathbf n}(\tilde f)\! \cdot\!
\mathrm{g}^{ \mathbf{ml} } \!\cdot\! \mathrm{g}^{ \mathbf{ik} }\right )
\cdot
C_{\mathbf {mn}}^{\mathbf {j}}(\tilde f)
\tilde \nabla_{\mathbf i} \tilde  \nabla_{\mathbf j},\label{bs}
\end{gather}
where $Sym$ denotes symmetrization.

Recall that the space $C^{\infty}(T^*M)$ is naturally a module over the group
$\mathrm{Diff}(M).$ The action is as follows: take $f\in \mathrm{Diff}(M)$ 
and $Q\in C^{\infty}(T^*M),$ then
\begin{equation}
\label{cam}
f^{*}\cdot Q:=Q \circ  \tilde f^{-1}.
\end{equation}
By differentiating the action above, one obtains the action of the Lie algebra
$\mathrm{Vect}(M).$

The action (\ref{cam}) induces an action of the group $\mathrm{Diff}(M)$ on
${\mathcal D} (C^{\infty}(T^*M))$ as follows: take
$f\in \mathrm{Diff}(M)$ and $T \in {\mathcal D} (C^{\infty}(T^*M)),$ then
\[
f^* T:= f^{*}\circ  T \circ {f^{-1}}^{*}.
\]
In this paper we shall study the cohomology arising in this context.

We have the following Theorem
\begin{theorem}
The map
\[
f \mapsto {\mathcal L} (f^{-1}),
\]
defines a non-trivial 1-cocycle on $\mathrm{Diff}(M)$ with values in
${\mathcal D}^3 (C^{\infty} (T^* M)).$
\end{theorem}

\begin{proof} To prove that the operator (\ref{bs}) is 1-cocycle one has to
verify the 1-cocycle condition
\begin{equation}
\label{moul}
{\mathcal L} (f\circ h)={\tilde h}^* {\mathcal L}(f)+ {\mathcal L}(h),
\end{equation}
for all $f,h\in \mathrm{Diff}(M).$

To verify (\ref{moul}) we use the formulae
\begin{gather}
Sym_{\mathbf i,\mathbf m,\mathbf n}
\left (C_{\mathbf k \mathbf l}^{\mathbf i}(\tilde f)
\!\cdot\! \mathrm{g}^{\mathbf{mk} }\! \cdot \mathrm{g}^{ \mathbf{nl} }\right )\!
\cdot\! C_{\mathbf m \mathbf n}^{\mathbf j}(\tilde f)
\tilde  \nabla_{\mathbf i} \tilde  \nabla_{\mathbf j}
({\tilde h}^{{-1}^*} Q)\nonumber\\
\qquad {}= Sym_{\mathbf i,\mathbf m,\mathbf n}
\left (C_{\mathbf k \mathbf l}^{\mathbf i}(\tilde f)
\!\cdot\! \mathrm{g}^{\mathbf{mk} }\! \cdot \mathrm{g}^{ \mathbf{nl} }\right )
\!\cdot\! C_{\mathbf m \mathbf n}^{\mathbf j}(\tilde f)
{\tilde h}^{{-1}^*} \tilde \nabla_{\mathbf i} \tilde
\nabla_{\mathbf j} Q,\nonumber\\
Sym_{\mathbf j,\mathbf i,\mathbf k} \left (C_{\mathbf m \mathbf l}^{\mathbf j}
(\tilde f) \!\cdot\! \mathrm{g}^{\mathbf{im} }\! \cdot \mathrm{g}^{ \mathbf{kl} }\right )
\tilde \nabla_{\mathbf i} \tilde  \nabla_{\mathbf j}
\tilde \nabla_{\mathbf k}({\tilde h}^{{-1}^*} Q)\nonumber\\
\qquad {}=Sym_{\mathbf j,\mathbf i,\mathbf k}
\left (C_{\mathbf m \mathbf l}^{\mathbf j}(\tilde f)
\! \cdot\!\mathrm{g}^{\mathbf{im} }\! \cdot \mathrm{g}^{ \mathbf{kl} }\right )\!
\left ( 3C_{\mathbf{ij}}^{\mathbf t}
({\tilde h}^{{-1}})
 {\tilde h}^{{-1}^*} \tilde \nabla_{\mathbf t} \tilde
\nabla_{\mathbf k} Q
+
 {\tilde h}^{{-1}^*} \tilde \nabla_{\mathbf i} \tilde
\nabla_{\mathbf j}  \tilde \nabla_{\mathbf k} Q\right )\!, \label{fo}
\end{gather}
for all $f, h \in \mathrm{Diff}(M)$ and for all $Q\in C^{\infty} (T^*M).$
\newpage

Let us prove that this 1-cocycle is not trivial. Suppose that there exists an
operator $A$ of degree three such that
\[
{\mathcal L}(f)={\tilde f}^* A-A.
\]
The principal symbol of the operator $A$ transforms under coordinates change
as a contravariant tensor field of degree three. Thus, it depends only on the
first jet of the symplectomorphism $\tilde f,$ whereas the principal symbol of
the operator (\ref{bs}) depends on the second jet of the symplectomorphism
$\tilde f,$ a contradiction.
\end{proof}

\begin{remark}
The fact that the first-order coefficients in the formulae (\ref{fo}) are
trivial is due to particular properties of the connection~(\ref{lift}).
\end{remark}

\subsection{Expression in local coordinates}

Denote by $(x^i,\xi_i), i=1,\ldots,n,$ the local coordinates on $T^*M.$ In
these coordinates the components of the symplectic form $\mathrm{g}$ are
\[
\mathrm{g}_{ij}= \mathrm{g}_{\bar i\bar j}=0,\qquad  \mathrm{g}_{i\bar j}= \delta_{ij}, \qquad
\mathrm{g}_{i\bar j}=-\mathrm{g}_{\bar ji},
\]
where $i,j=1,\ldots,n.$

For any diffeomorphism $f(x)=\left(f^1(x),\ldots,f^n(x)\right)$ on $M,$ its 
symplectic lift to $T^*M$ is written as $\tilde f(x,\xi)=\left( f^1(x),\ldots, 
f^n(x),\frac{\partial x^i}{\partial f^1} \xi_i, \ldots,
\frac{\partial x^i}{\partial f^n} \xi_i\right)$.

An easy computation shows that the operator (\ref{bs}) (to some factor)
has the form
\begin{gather}
3 C^i_{jk}(\tilde f)
 \frac{\partial^2}{\partial \xi^j\partial \xi^k } \frac{\partial}{\partial x^i}
 +3 \left (2
C^m_{ik}(\tilde f)\Gamma^k_{jm}+ C^m_{ki}(\tilde f)C^k_{mj}(\tilde f)
\right) \frac{\partial^2}{\partial \xi^i \partial \xi^j}\nonumber\\
\qquad {}+ C^{\bar i}_{jk}(\tilde f) \frac{\partial^3}{\partial \xi^i
\partial \xi^j\partial \xi^k},\label{exp1}
\end{gather}
where $C(\tilde f)$ is the tensor (\ref{masri}).

In the particular case when the connection $\Gamma$ is Euclidean
(i.e. $\Gamma \equiv 0$), the operator above takes the form
\begin{gather}
3\frac{\partial^2 f^l}{\partial x^j \partial x^k}
\frac{\partial x^i}{\partial f^l}
\frac{\partial^2}{\partial \xi^j\partial \xi^k } \frac{\partial}{\partial x^i}+
3\frac{\partial^2 f^k}{\partial x^q \partial x^i}
\frac{\partial x^m}{\partial f^k}
\frac{\partial^2 f^l}{\partial x^j \partial x^m}
\frac{\partial x^q}{\partial f^l}
\frac{\partial^2}{\partial \xi^i \partial \xi^j}\nonumber\\
\qquad {}+\frac{\partial x^m}{\partial f^q}\xi_m
\left (3\frac{\partial x^l}{\partial f^p}
\frac{\partial^2 f^p}{\partial x^i \partial x^j}
\frac{\partial^2 f^q}{\partial x^l \partial x^k}
-\frac{\partial^3 f^q}{\partial x^i \partial x^j \partial x^k}\right )
\displaystyle \frac{\partial^3}{\partial \xi^i \partial \xi^j\partial \xi^k}.
\label{exp2}
\end{gather}

\section{Cohomology of the group of diffeomorphisms}

Let $M$ be an oriented manifold. Denote by $\mathrm{Diff}_+ (M)$ the group of
diffeomorphisms of $M$ that preserve the orientation on $M.$

\subsection{Space of linear differential operators and space of symbols}

Let ${\mathcal D}({\mathcal F}_\lambda)$ be the space of linear differential 
operators acting on the space of tensor densities of degree $\lambda$ on $M.$ 
This space admits a structure of module over the group $\mathrm{Diff}(M)$ as 
follows:  take 
$f\in \mathrm{Diff}(M)$ and $A\in {\mathcal D} ({\mathcal F}_{\lambda}),$ then
\begin{equation}
\label{op}
f^* A:= f^*\circ A \circ {f^{-1}}^{*},
\end{equation}
where $f^*$ is the natural action of a diffeomorphism on $\lambda$-densities.

By differentiating this action, one obtains the action of $\mathrm{Vect}(M)$ on
${\mathcal D}({\mathcal F}_{\lambda})$.

Consider the space ${\mathcal S}^k$ of symmetric contravariant tensor
fields on $M$ of degree $k.$ This space is naturally isomorphic
to the space of functions on $T^* M$ that are polynomials of degree $k$
on the fibers. We shall identify these two spaces throughout this paper.

One can define a $\mathrm{Diff}(M)$-module structure on ${\mathcal S}^k$ by the
formula (\ref{cam}). We have then a graduation of $\mathrm{Diff}(M)$-module
\[
{\mathcal S}=\bigoplus_{k\geq 0}{\mathcal S}^k.
\]
The module ${\mathcal S}$ and the modules ${\mathcal D}({\mathcal F}_{\lambda})$
are not isomorphic with respect to the action of $\mathrm{Diff}(M)$ given 
above (cf. \cite{boh, do, lo,lo2}).

\subsection{Cohomology of $\boldsymbol{\mathrm{Vect} (M)}$
and cohomology of $\boldsymbol{\mathrm{Diff}(M)}$}
\label{alg}

The problem of the infinitesimal deformation of the module $\mathcal S$ with
respect to the Lie algebra $\mathrm{Vect}(M)$ is related to the cohomology 
group $\mathrm{H}^1(\mathrm{Vect}(M); \mathrm{End} ({\mathcal S})).$ Consider 
now the space of differential operators ${\mathcal D} (\mathcal S)$ as 
a submodule of $\mathrm{End} ({\mathcal S}).$
This module can be decomposed as a $\mathrm{Vect} (M)$-module into
$\bigoplus_{k,l} {\mathcal D} ({\mathcal S }^k, {\mathcal S}^l)$. The 
following result is proved in~\cite{lo2}.
\begin{equation}
\label{cal}
\mathrm{H}^1(\mathrm{Vect}(M); {\mathcal D}({\mathcal S}^k,{\mathcal S}^m))=
\left\{
\begin{array}{ll}
{\mathbb R},& \mbox{if} \quad k-m=2,\\
{\mathbb R},& \mbox{if} \quad k-m=1,\ m\not=0,\\
{\mathbb R}\oplus {\mathrm H}^1_{\mathrm{DR}}(M),& \mbox{if} \quad k-m=0,\\
0,&\hbox{otherwise.}
\end{array}
\right.
\end{equation}
Let us give explicit formulae for 1-cocycles that generate the cohomology
group (\ref{cal}).

{\bf For  $\boldsymbol{k-m=0}$.}
Denote by $\mathrm{Div}$ the divergence operator associated with the volume 
form on $M.$ The map $X\mapsto a \,\mathrm{Div}(X)+ i_X\omega,$
where $a$ is a  constant and $\omega$ is a closed 1-form on $M,$ defines a
1-cocycle on $\mathrm{Vect}(M)$ with values in $C^{\infty}(M)$ (cf.~\cite{f})

The 1-cocycle generating the cohomology group (\ref{cal}), for
$k-m=0,$ is a multiplication by $k$-tensor fields with the cocycle above.

{\bf For $\boldsymbol{k-m=1,}$ and $\boldsymbol{m\not=0}$.}
The map $X\mapsto L_X \nabla,$ where
$ L_X \nabla$ is the Lie derivative of the connection $\nabla$ (see~\cite{lo2}
for more details), defines a 1-cocycle on $\mathrm{Vect}(M)$ with values in the space
of tensor fields on $M$ of type $(2,1).$

The 1-cocycle generating the cohomology group~(\ref{cal}),
for $k-m=1,$ is obtained by contracting $k$-tensor fields with the above
(2,1)-tensor fields.

{\bf For  $\boldsymbol{k\!-\!m=2}$.}
It is well known that the cohomology group
$\mathrm{H}^2 (C^{\infty}(T^*M); C^{\infty}(T^*M))$ is isomorphic to
${\mathbb R} \oplus \mathrm H^1_{\mathrm{DR}}(M)$
(see, e.g.~\cite{or} and the references therein). The component 
${\mathbb R}$ is generated by the so-called Vey cocycle~\cite{v}, usually 
denoted by $S_{\Gamma}^3$ (see \cite{ro} for an explicit construction which 
uses a connection). The Vey cocycle is important in the theory of deformation
quantization; it appears in the third-order term in the Moyal product
(see, e.g.~\cite{or}).

Consider now the natural embedding of $\mathrm{Vect}(M)$ in 
$C^{\infty} (T^*M)$ given by $X\mapsto X^i\xi_i.$ The map
$X\mapsto S_{\tilde \Gamma}^3 (X^i\xi_i,\cdot),$ where $\tilde \Gamma$ is the
connection (\ref{lift}), turns out to be a 1-cocycle on $\mathrm{Vect}(M)$ 
with values in ${\mathcal D}({\mathcal S}^k,{\mathcal S}^{k-2})$ and 
generates the cohomology group (\ref{cal}) (see \cite{lo2}).
\begin{remark}
In the one dimensional case, the cohomology group (\ref{cal}) is
calculated in \cite{bo} (see \cite{bg} for the complex case).
\end{remark}

To study the cohomology of the group of diffeomorphisms, we deal with
differential cohomology ``Van Est Cohomology''. This means we consider only
differential cochains (see~\cite{f}).

We have the following:
\begin{theorem}
\label{salam}
Let $M:=S^{n}$ be the $n$-dimensional sphere. For $n=2,3,$ the
first-cohomology group
\begin{equation}
\label{main}
\mathrm{H}^1(\mathrm{Diff}_{+}(S^n); {\mathcal D}({\mathcal S}^k,{\mathcal S}^m))=
\left\{
\begin{array}{ll}
{\mathbb R},&\mbox{if} \quad k-m=0,\\
{\mathbb R},&\mbox{if} \quad k-m=2,\\
{\mathbb R},&\mbox{if} \quad k-m=1,m\not=0,\\
0,&\mbox{otherwise.}
\end{array}
\right.
\end{equation}
\end{theorem}
We start by giving explicit formulae for 1-cocycles generating the
cohomology group (\ref{main}).

{\bf For  $\boldsymbol{k-m=0}$.}
Any diffeomorphism $f\in \mathrm{Diff}_{+}(S^n)$ preserves the volume form on 
$S^n$ up to some factor. The logarithm function of this factor defines a 
1-cocycle on $\mathrm{Diff}(S^n),$ say $c(f),$ with values in 
$C^{\infty}(S^n).$ 

The 1-cocycle generating the component ${\mathbb R}$ of the cohomology group
(\ref{main}) is given as a multiplication operator of any $k$-tensor field by
the 1-cocycle $c(f).$

{\bf For  $\boldsymbol{k-m=1}$, $\boldsymbol{m\not=0}$.} The
difference $\ell(f):=f^* \Gamma -\Gamma,$ where $\Gamma$ is a
connection on~$S^n,$ is a well-defined tensor field of type
$(2,1).$ The map $f\mapsto \ell(f^{-1})$ defines a 1-cocycle on
$\mathrm{Diff}_+(S^n)$ with values in the space of tensor fields
of type $(2,1).$

The 1-cocycle generating the cohomology group (\ref{main}) is obtained by
contracting every $k$-tensor field with the above 1-cocycle.

{\bf For  $\boldsymbol{k-m=2}$.}
First we have the following
\begin{proposition}
\label{doc}
For any $f\in \mathrm{Diff}_+(S^n),$ the restriction
${\mathcal L}(\tilde f)_{|_{{\mathcal S}^k}},$ where ${\mathcal L}$ is the
opera\-tor~(\ref{bs}), defines a map from ${\mathcal S}^k$ to 
${\mathcal S}^{k-2}.$ 
\end{proposition}

The proof is immediate from the formula (\ref{exp1}) and the
properties of the connec\-tion~(\ref{lift}).

It is easy to see that the map
$f\mapsto {\mathcal L}({\tilde f^{-1}})_{|_{{\mathcal S}^k}}$
defines a 1-cocycle on $\mathrm{Diff}_+(S^n)$ with values in
${\mathcal D}({\mathcal S}^k, {\mathcal S}^{k-2}),$ and generates the 
cohomology group~(\ref{main}).

\begin{remark} It is worth noticing that the 1-cocycles above can
be defined on any oriented manifold.
\end{remark}

Now we are in position to prove Theorem~\ref{salam}.

It is well-known that the maximal compact group of ``rotations'' of $S^n$,
$SO(n+1),$ is a deformation retract of the group $\mathrm{Diff}_+ (S^n),$
for $n=1,2,3$ (see~\cite{th}). Since the space $\mathrm{Diff}_+ (S^n)/SO(n+1)$ 
is acyclic, the Van Est cohomology of the Lie
group $\mathrm{Diff}_+(S^n)$ can be computed using the isomorphism
(see, e.g., \cite[p.~298]{f})
\[
H^1(\mathrm{Diff}_{+}(S^n); {\mathcal D}({\mathcal S}^k,{\mathcal S}^m))\simeq
H^1(\mathrm{Vect} (S^n), SO(n+1);{\mathcal D}({\mathcal S}^k,{\mathcal S}^m)).
\]
Thus $H^1(\mathrm{Diff}_{+}(S^n); {\mathcal D}({\mathcal S}^k,{\mathcal S}^m))=0$
for $k-m\neq 0, 1, 2,$ and for $(k,m)=(1,0).$

First, observe that the De Rham classes in the cohomology group (\ref{cal})
turn out to be trivial since $H^1_{\mathrm {DR}}(S^n)=0.$ Suppose
that there are two 1-cocycles representing cohomology classes
in the cohomology group~(\ref{main}). The isomorphism above shows that these
two 1-co\-cycles induce two non-cohomologue classes in the cohomology
group (\ref{cal}), which is absurd. Theorem~\ref{salam} follows from explicit
constructions of the 1-cocycles above.\hfill\qed

\begin{remark} (i) Theorem \ref{salam} remains true for $n\geq 4$ if the
subgroup $SO(n+1)$ is a deformation retract of the group 
$\mathrm{Diff}_+(S^n)$ for all $n.$ As far as I know, this problem is not 
solved yet.

(ii) It is possible to integrate the De Rham class to the group 
$\mathrm{Diff}_+(M),$ when $M$ is an (oriented) simply connected manifold, if
one considers non-differential cochains. To define these 1-cocycles,
let $\alpha$ be a closed 1-form on $M.$ The function
\begin{equation}
\label{derham}
c^{\alpha}(f)(x):=\int_{\gamma (x,f(x))} \alpha,
\end{equation}
where $\gamma (x,f(x))$ is a path joining the point $x$ to the point $f(x),$
is well-defined on $M$ since $\alpha$ is closed and $M$ is simply connected.
It is easy to show that the map
$f\mapsto c^{\alpha}(f^{-1})(x)$ defines a 1-cocycle on $\mathrm{Diff}_{+}(M)$ with values in $C^{\infty}(M).$ The contraction of any tensor by the 1-cocycle
above defines a 1-cocycle that integrate, for $k-m=0$, the De Rham class
in the cohomology group (\ref{cal}).
\end{remark}

\section{Discussion}
\label{las}

{\bf 1.} In the one dimensional case $M:=S^1,$ the
analogue of Theorem~(\ref{main}) was given in \cite{bo}: the first
cohomology group
\begin{equation}
\label{dim}
\mathrm{H}^1 (\mathrm{Diff}_+(S^1),
\mathrm{PSL}_2({\mathbb R}); {\mathcal D} ({\mathcal F}_\lambda, {\mathcal F}_\mu))=
\left\{
\begin{array}{ll}
{\mathbb R},&\mbox{if} \quad \mu-\lambda=2, 3, 4,\ (\lambda \mbox{ generic})\\
{\mathbb R},&\mbox{if} \quad (\lambda,\mu)=(-4,1),(0,5)\\
0,&\hbox{otherwise.}
\end{array}
\right.
\end{equation}
The 1-cocycles generating the cohomology group~(\ref{dim}) are
differential operators whose coefficients are a linear combination
by the so-called Schwarzian derivative and its derivatives
(see~\cite{boh} for explicit formulae).
 
From this point of view, Ovsienko and the author~\cite{boh}
(see also~\cite{bh}) generalize the Schwarzian derivative on $S^n$
endowed with a projective structure (i.e.~coordinates change
are projective transformations): the multi-dimensional Schwarzian
derivative is a 1-cocycle on the group $\mathrm{Diff}(S^n)$ with values in
${\mathcal D}({\mathcal S}^2, {\mathcal S}^0)$,
vanishing on the subgroup $\mathrm{PSL}_{n+1}({\mathbb R}).$ It is 
interesting to ask whether there is a relation between the 1-cocycle 
(\ref{exp2}) and the multi-dimensional Schwarzian derivative. Observe that 
the 1-cocycle (\ref{exp2}) does not vanish on the subgroup 
$\mathrm{PSL}_{n+1}({\mathbb R})$.

Another approach for the {\it conformal\/} Schwarzian derivative
on $S^n$ was introduced in \cite{bs}. This derivative is a
1-cocyle on the group $\mathrm{Diff}(S^n)$ with values in
${\mathcal D}({\mathcal S}^2, {\mathcal S}^0)$, vanishing on the
subgroup $\mathrm{O}(p+1,q+1)$ (group of all conformal
diffeomorphisms with respect to a given metric on $S^n$, where
$p+q=n$). It is interesting to ask whether there is a relation
between the 1-cocycle~(\ref{exp2}) and the {\it conformal}
Schwarzian derivative introduced in~\cite{bs}.

{\bf 2.} It is interesting to generalize the 1-cocycle (\ref{bs}) to any
arbitrary symplectic manifold. In this case, the 1-cocycle (\ref{bs}) can be
interpreted as a cocycle integrating on one argument the Vey 2-cocycle
$S^3_{\Gamma}$ (see Section~\ref{alg}). The 1-cocycle (\ref{bs}) integrates
the Vey cocycle only on the cotangent bundle $T^*M.$

Another 2-cocycle was constructed by Tabachnikov \cite{tab} integrating
the Vey 2-cocycle. It is also interesting to compare this 2-cocycle with the
operator~(\ref{bs}).

\subsection*{Acknowledgments}

I am grateful to V~Ovsienko for his
constant support and for helpful remarks concerning this paper. I am also
grateful to H~Gargoubi, B~Kolev, Y~Maeda, P~Lecomte and  H~Moriyoshi
for numerous fruitful discussions.
Research supported by the Japan Society
for the Promotion of Science.

 \label{Bouarroudj-lastpage}

\end{document}